\documentclass{elsart}

\usepackage{amsmath}
\usepackage{amssymb}
\journal{Adv. Appl. Math.} 

\newenvironment{GT}{\textbf{Gordon's Partition Theorem} \em  }{}
\newenvironment{AG}{\textbf{The Andrews-Gordon Theorem} \em }{}

\newcommand{\hgs}[6]{ {}_{#1}\phi_{#2} \left[ \genfrac{}{}{0pt}{}{#3}{#4} ; {#5},{#6} \right]}
\newcommand{\vwp}[6]{ {}_{#1}W_{#2} \left( {#3};{#4};{#5},{#6} \right)}
\newcommand{\binomial}[2]{ \genfrac{(}{)}{0pt}{}{#1}{#2} }

\numberwithin{equation}{section}
\numberwithin{thm}{section}

\begin{document}
\begin{frontmatter}
\title{On Simplifications of Certain $q$-Multisums}
\author{Andrew V. Sills}
\address{Department of Mathematics, Rutgers University,
         Hill Center, Busch Campus, Piscataway, NJ 08854-8019 USA}
\ead{asills@math.rutgers.edu}
\ead[url]{www.math.rutgers.edu/\~{}asills}
\date{December 22, 2005; accepted June 2, 2006}

\begin{abstract}
  Some examples of naturally arising 
multisum $q$-series which turn out to
have representations as 
fermionic 
single sums are presented. The resulting identities are proved using
transformation formulas from the theory of basic hypergeometric series.
\end{abstract}


\end{frontmatter}

\section{Introduction}

It is certainly useful to be able to transform a naturally arising
$q$-multisum into a single-fold sum. 
Many such identities are well known, e.g. the ``$a$-generalized Andrews-Gordon theorem,"
(see Eq.~\eqref{aGG} below), and Bressoud's generalization to even moduli~\cite{DMB1,DMB2},
Krattenthaler and Rosengren's $q$-analog~\cite{CK} of an identity of 
Gelfand, Graev, and Retakh~\cite{GGR}, the many identities of
the physicists Berkovich, McCoy, Orrick, Pearce, Schilling, and 
Warnaar~\cite{phys1,phys2,phys3,phys4,phys5,phys6}, to name a few.
In these papers, the multisums generally fall into
a category of what physicists call ``fermionic" representations, while
the single sum representations are ``bosonic."

 A different category of $q$-identities, often called
``fermionic reduction formulas," first appears in Andrews'
1981 paper on 
multiple $q$-Series identities~\cite{GEA:Multiq}, where 
it is shown how to simplify certain fermionic $q$-multisums via
two ``amalgamation lemmas"~\cite[pp. 19, 20; Lemmas 1, 2]{GEA:Multiq}. 
Another important result in this genre was conjectured by
Melzer~\cite{EM}, and proved by Bressoud, Ismail, and 
Stanton using Bailey lemma techniques ~\cite[Thms. 5.1 and 5.2]{BIS}.  
The Melzer conjecture
was subsequently reproved by Warnaar~\cite[Thm. 4.4]{SOW2} using
other methods.  In~\cite{SOW1}, Warnaar provides another 
interesting fermionic reduction formula. 
The results of this paper fall into the category of fermionic
reduction formulas.

%\begin{NoteTerm}
  The terms ``fermionic" and ``bosonic" arise in statistical physics.
  For proper definitions of ``fermionic" and ``bosonic," one should 
consult an appropriate paper written by a physicist; 
see e.g. ~\cite[p. 165 ff.]{BM:Century}.
However, since $q$-series identities arising in physics also occur
in combinatorics and mathematical analysis, it seems reasonable,
by analogy,
to attach the term ``fermionic" to any 
set of integer partitions with difference conditions and to its
associated generating function, and the term ``bosonic" to
a set of integer partitions associated with congruence conditions,
and to its associated generating function.  This is the sense in
which I use the terms herein. 
%\end{NoteTerm}

 Before proceeding to the identities, it would be useful to point out the
combinatorial context which paved the way to their discovery. 

 A \emph{partition} $\pi$ of an integer $n$ is a nonincreasing sequence $(\pi_1, \pi_2, \pi_3, \dots )$ of nonnegative integers
such that $\sum_{j=1}^\infty \pi_j=n$.  Each nonzero $\pi_i$ is called a \emph{part} of $\pi$.  The number of times 
$j$ appears in $\pi$ is called the \emph{multiplicity} of $j$ in $\pi$ and is denoted $m_j(\pi)$.

Recall Gordon's combinatorial generalization of the Rogers-Ramanujan identities~\cite{BG}:

\begin{GT} Let $A_{k,i} (n)$ denote the number of partitions of the integer $n$ into parts  $\not\equiv 0,\pm i \pmod{2k+1}$.
Let $B_{k,i}(n)$ denote the number of partitions 
$\pi$ of $n$ such that $m_1(\pi)<i$, and for any positive integer $j$,
$m_j(\pi)+m_{j+1}(\pi)<k$.  Then for $1\leqq i\leqq k$, $A_{k,i}(n) = B_{k,i}(n)$.
\end{GT}

It is well known~\cite[p. 111, Theorem 7.8]{GEA:TOP} that Andrews' 
generalization of the Rogers-Ramanujan
identities for odd moduli~\cite[p. 4082, Theorem 1]{GEA:OddModuli} is a $q$-series counterpart to Gordon's theorem, and as such is often called
the ``Andrews-Gordon theorem." 

\begin{AG} For $1\leqq i \leqq k$, $k\geqq 1$, 
\begin{equation}\label{AndGor}
\sum_{n_1, n_2, \dots, n_{k-1}\geqq 0} 
\frac{q^{N_1^2 + N_2^2+ \cdots+ N_{k-1}^2 + N_i + N_{i+1} + \cdots + N_{k-1}}}
{(q)_{n_1} (q)_{n_2} \cdots (q)_{n_{k-1}} } 
= \frac{(q^i, q^{2k+1-i}, q^{2k+1}; q^{2k+1})_\infty}{(q)_\infty},
\end{equation}
where $N_j = n_j + n_{j+1} + \cdots + n_{k-1}$, and 
\[ (a)_n = (a;q)_n = (1-a)(1-aq)(1-aq^2)\cdots(1-aq^{n-1}), \]
\[ (a)_\infty = (a;q)_\infty = (1-a)(1-aq)(1-aq^2) \cdots, \]
\[ (a_1, a_2, \dots, a_r; q)_n = (a_1)_n (a_2)_n \cdots (a_r)_n. \]
\end{AG}

 Furthermore, a standard refinement of the $B_{k,i}(n)$ of Gordon's partition theorem counts $\mathcal{B}_{k,i} (m,n)$, the number of
partitions of $n$ counted by $B_{k,i}(n)$ which have exactly $m$ parts.  The analogous refinement of 
\eqref{AndGor} 
is~\cite[p. 112, Eq. (7.3.8)]{GEA:TOP}
\begin{multline}
\sum_{n=0}^\infty \sum_{m=0}^n \mathcal B_{k,i} (m,n) a^m q^n \\= 
\sum_{n_1, n_2, \dots, n_{k-1}\geqq 0} 
\frac{a^{N_1+N_2+\cdots+N_{k-1}}
q^{N_1^2 + N_2^2+ \cdots+ N_{k-1}^2 + N_i + N_{i+1} + \cdots + N_{k-1}}}
{(q)_{n_1} (q)_{n_2} \cdots (q)_{n_{k-1}} }\\
=\frac{1}{(aq)_\infty} \sum_{n\geqq 0} 
\frac{ (-1)^n a^{kn} q^{ (2k+1)n(n+1)/2-in } (1-a^i q^{(2n+1)i}) (aq)_n }{(q)_n} \label{aGG}
\end{multline}

In a recent paper~\cite{AVS:rrt}, I showed that certain $q$-series were related to dilated versions of special cases of Gordon's partition theorem.  
Accordingly,
the following identities between multisum and single sum $q$-series, although not explicitly stated in~\cite{AVS:rrt}, follow immediately from the
results therein:

\begin{gather}
  \sum_{r=0}^\infty \sum_{s=0}^\infty \frac{a^{r+2s} 
  q^{2(r+s)^2 + 2s^2}}{(q^2;q^2)_{r} (q^2;q^2)_{s}}
 = (aq;q^2)_\infty \sum_{j=0}^\infty \frac{ a^j q^{j^2} }{ (aq;q^2)_{j} (q)_j}
\label{d2k3i3}\\
  \sum_{r=0}^\infty \sum_{s=0}^\infty \sum_{t=0}^\infty
     \frac{a^{r+2s+3t} q^{3(r+s+t)^2 + 3(s+t)^2 + 3t^2}}{(q^3;q^3)_r (q^3;q^3)_s (q^3;q^3)_t}
     =\frac{(aq)_\infty}{(aq^3;q^3)_\infty}\sum_{j=0}^\infty \frac{a^j q^{j^2} (a;q^3)_j}{(a)_{2j} (q)_j}
     \label{d3k4i4}.
\end{gather}

 Note that the series on the left hand sides of \eqref{d2k3i3} and \eqref{d3k4i4} are special cases of 
the Andrews-Gordon theorem (with $q\to q^2$, $k=i=3$; and $q\to q^3$, $k=i=4$ respectively), while
the right hand sides arise in Bailey's two-variable generalizations~\cite[pp. 6--7]{WNB1948}
of Rogers' first mod 14 identity~\cite[p. 341, Ex. 2]{LJR1894} 
and Dyson's first mod 27 identity~\cite[p. 433, Eq. (B4)]{WNB1947} respectively.

 Note that the Andrews-Gordon theorem, and
indeed all of the identities presented herein, may be regarded as
identities of analytic functions and are thus subject to convergence conditions.
However, since the underlying motivation is combinatorial (and so the series may 
be regarded
as generating functions), the convergence
conditions 
will not be explicitly mentioned.

%\begin{rem} 
Although proved in \cite{AVS:rrt} with the aid of systems of
$q$-difference equations, 
it is not at all obvious \emph{why} the right hand side 
of \eqref{d2k3i3} enumerates the partitions from the $k=i=3$ case of Gordon's partition theorem (dilated by a factor of 2)
and the right hand side of~\eqref{d3k4i4} enumerates the parititons from the $k=i=4$ case of Gordon's partition
theorem (dilated by a factor of 3), but once this fact is established, 
their equality with their respective left hand sides
is immediate thanks to the Andrews-Gordon theorem.
%\end{rem}

  The purpose of this note is to present $q$-hypergeometric 
proofs of \eqref{d2k3i3} and \eqref{d3k4i4} in order
to gain an understanding of these identities from the standpoint of basic hypergeometric series.
The identities \eqref{d2k3i3} and \eqref{d3k4i4} will be derived as corollaries of the more general identities

\begin{multline} \label{gen14}
    \sum_{j=0}^n \sum_{h=0}^{j} \frac{(q^{-2n})_{2j} (yq;q^2)_j (aq/xy;q^2)_j
      (q^{-2j};q^2)_h (xq;q^2)_h (y;q^2)_h  }
    {(q^2;q^2)_j (aq/x;q^2)_j (aq^2/y;q^2)_j (yq^{2-4n}/a;q^2)_j (q^2;q^2)_h (aq^2/b;q^2)_h
     }\\
     \times \frac{(aq^2/bx;q^2)_h } { (aq^2/x;q^2)_h (xyq^{1-2j};q^2)_h } q^{2h+2j}\\
   = \frac{(aq;q^2)_{2n} (aq/xy)_{2n}}{(aq/x)_{2n} (aq^2/y;q^2)_{2n-1} (1-aq^{2n}/y)} 
  \sum_{j=0}^{2n} \frac{ (x)_j (y)_j (aq/b;q^2)_j (q^{-2n})_j q^j}
 {(q)_j (aq;q^2)_j (aq/b)_j (xyq^{-2n}/a)_j}
 \end{multline}
 and
 \begin{multline}\label{gen27}
 \sum_{j=0}^{n} \sum_{k=0}^j \sum_{m=0}^k
   \frac{ (q^{-3n})_{3j} (a/y^2;q^3)_j (q^{-3j};q^3)_k (yq;q^3)_k (yq^2;q^3)_k 
     }
   { (q^3;q^3)_j (q^3;q^3)_k (q^3;q^3)_m (aq^2/y; q^3)_j (aq/y;q^3)_j (aq/x;q^3)_k }
   \\
  \quad \times \frac{(aq/xy;q^3)_k (q^{3k};q^3)_m (xq^2;q^3)_m (y;q^3)_m (aq^2/x^2;q^3)_m q^{3j+3k+3m}} { (aq^2/x;q^3)_m (aq^3/x;q^3)_m  (q^{3-9n};q^3)_j (y^3 q^{3-3j}/a;q^3)_k (xy q^{2-3j}/a;q^3)_m
  (aq^3/y;q^3)_k}\\
   = \frac{  (aq)_{3n} (aq/xy)_{3n}  (aq^{3n+1};q^3)_n (aq^{3n+2};q^3)_n} 
     {(aq/x)_{3n} (aq/y)_{3n}  (aq^3;q^3)_n (aq^{6n};q^3)_n  }
     \sum_{j=0}^{3n} \frac{  (a;q^3)_j (x)_j (y)_j (q^{-3n})_j q^j} { (q)_j  (a;q^2)_{2j}  (xyq^{-3n}/a)_j }.
 \end{multline}

 In \S\ref{Proofs}, proofs of~\eqref{d2k3i3}--\eqref{gen27} are presented.  
These proofs suggest additional results, presented in \S\ref{AddResults}.  
Finally, related open questions are presented in \S\ref{Disc}. 

\section{$q$-Hypergeometric proofs of~\eqref{d2k3i3}--\eqref{gen27}}\label{Proofs}
The basic hypergeometric series is defined, as in~\cite{GR}, by
\begin{multline} \label{bhsdef} \hgs{r}{s}{a_1, a_2, \dots, a_r}{b_1, b_2, \dots, b_s}{q}{z} \\ :=
\sum_{j=0}^\infty \frac{(a_1;q)_j (a_2;q)_j \cdots (a_r;q)_j}
{(q;q)_j (b_1;q)_j (b_2;q)_j \cdots (b_s;q)_j}z^j \left[ (-1)^j q^{j(j-1)/2} 
\right]^{1+s-r},
\end{multline}

Note that a basic hypergeometric series~\eqref{bhsdef} is called 
\emph{well-poised} if 
\mbox{$s=r-1$} and  $a_1 q= a_2 b_2 = \cdots=a_{r} = b_{r-1}$ and 
\emph{very-well-poised} if,
in addition, $a_2 = q\sqrt{a_1}$ and $a_3 = -q\sqrt{a_1}$.  
Very-well-poised series are central
to the study and turn out to be the common link between the multisum and single sum representations
in this paper.  It will be convenient to employ the following condensed notation for very-well-poised
basic series:
\begin{multline*}
\vwp{r+1}{r}{a}{a_4, a_5, \dots, a_{r+1}}{q}{z} \\ :=
     \hgs{r+1}{r}{a, q\sqrt{a},-q\sqrt{a},a_4, a_5,\dots, a_{r+1}}
     {\sqrt{a},-\sqrt{a},a q/a_4, a q/a_5, \dots, aq/a_{r+1}}{q}{z}.
 \end{multline*}
  
   In~\cite{GEA:ProbPros}, Andrews presents a very general series transformation formula 
whereby a very-well-poised 
${}_{2k+4}\phi_{2k+3}$ is transformed into a $(k-1)$-fold multisum
representation.  The $k=2$ case is equivalent to Watson's $q$-analogue of Whipple's theorem~\cite{GNW}.
The $k=3$ case of Andrews' transformation may be stated as
 \begin{multline} \label{GEAk3}
   \vwp{10}{9}{a}{b,c,d,e,f,g,q^{-n}}{q}{a^3 q^{n+3}/bcdefg} \\
   =\frac{(aq)_n (aq/fg)_n}{(aq/f)_n (aq/g)_n} 
    \sum_{j=0}^n \frac{(q^{-n})_j (f)_j (g)_j (aq/de)_j q^j}{(q)_j (aq/d)_j (aq/e)_j (fgq^{-n}/a)_j}\\
    \times
    \hgs{4}{3}{q^{-j},d,e,aq/bc}{aq/b, aq/c, deq^{-j}/a}{q}{q},
  \end{multline}
where, here and throughout, $n$ is a nonnegative integer.  
  While~\eqref{GEAk3} transforms a fairly general very-well-poised 
${}_{10}\phi_9$ to a double sum, 
there are a number of transformation formulas known which transform a somewhat more specialized
very-well-poised ${}_{10}\phi_9$ to a single sum.  
For instance, Verma and Jain~\cite[p. 232,  Eq. (1.4)]{VJ} found
\begin{multline}\label{VJ10-9-2}
 \vwp{10}{9}{a}{b,x,xq,y,yq,q^{1-n},q^{-n}}{q^2}{a^3 q^{2n+3}/bx^2y^2} \\
 = \frac{(aq)_n (aq/xy)_n}{(aq/x)_n (aq/y)_n} 
 \hgs{5}{4}{x,y,\sqrt{aq/b},-\sqrt{aq/b},q^{-n}}{\sqrt{aq},-\sqrt{aq},aq/b,xyq^{-n}/a}{q}{q}
\end{multline}
 Note: the $n\to\infty$ case of~\eqref{VJ10-9-2} is given by Bailey (in a somewhat disguised form) as
 \cite[p. 6, Eq. (6.3)]{WNB1948}.
 With~\eqref{GEAk3} and~\eqref{VJ10-9-2} in hand, 
it is now time to establish~\eqref{gen14}.
 
\begin{thm} Identity~\eqref{gen14} is valid.
\end{thm}
\begin{pf}
 The result follows from the observation that 
 \begin{equation}\label{vwplink14}
  \vwp{10}{9}{a}{b,x,xq,y,yq,q^{1-2n},q^{-2n}}{q^2}{a^3 q^{4n+3}/bx^2 y^2} 
\end{equation}
can be transformed via either~\eqref{GEAk3} or~\eqref{VJ10-9-2}.
Transforming~\eqref{vwplink14} via~\eqref{GEAk3} yields 
\begin{multline}\label{2.5}
\frac{(aq^2;q^2)_n (aq^{2n}/y;q^2)_n}{(aq/y;q^2)_n (aq^{1+2n};q^2)_n} 
    \sum_{j=0}^n \frac{(q^{-2n};q^2)_j (yq;q^2)_j (q^{1-2n};q^2)_j (aq/xy;q^2)_j}
    {(q^2;q^2)_j (aq/x;q^2)_j (aq^2/y;q^2)_j (yq^{2-4n}/a;q^2)_j}\\ \times
    \hgs{4}{3}{q^{-2j},xq,y,aq^2/bx}{aq^2/b, aq^2/x, xyq^{1-2j}/a}{q^2}{q^2},
\end{multline} while transforming~\eqref{vwplink14} via~\eqref{VJ10-9-2}
yields 
\begin{equation}\label{2.6}
\frac{(aq)_{2n} (aq/xy)_{2n}}{(aq/x)_{2n} (aq/y)_{2n}} 
 \hgs{5}{4}{x,y,\sqrt{aq/b},-\sqrt{aq/b},q^{-2n}}{\sqrt{aq},-\sqrt{aq},aq/b,xyq^{-2n}/a}{q}{q}.
\end{equation}
Thus~\eqref{2.5}$=$~\eqref{2.6}.
\qed\end{pf}

While it must be admitted that Identity~\eqref{gen14} is probably not the most 
beautiful of identities, it nonethless gives rise to elegant corollaries, which may
now be easily deduced.
\begin{cor} Identity~\eqref{d2k3i3} is valid.\end{cor}
\begin{pf} Let $b, x, y, n\to\infty$ in Eq.~\eqref{gen14}. \qed\end{pf}

Actually,~\eqref{d2k3i3} is just one of a set of three closely related
identities.  With~\eqref{d2k3i3} established, it is straightforward to deduce
its two partners:

\begin{cor}
\begin{gather}
 \sum_{r=0}^\infty \sum_{s=0}^\infty 
   \frac{a^{r+2s} q^{2(r+s)^2 + 2s^2 + 2r+4s}}{(q^2;q^2)_r (q^2;q^2)_s}
 = (aq;q^2)_\infty \sum_{j=0}^\infty \frac{a^j q^{j^2+2j}}{(aq;q^2)_{j+1} (q)_j}
  \label{d2k3i1}\\
  \sum_{r=0}^\infty \sum_{s=0}^\infty
   \frac{a^{r+2s} q^{2(r+s)^2 + 2s^2 + 2s}}{(q^2;q^2)_r (q^2;q^2)_s}
 = (aq;q^2)_\infty \sum_{j=0}^\infty \frac{a^j q^{j^2+j}}{(aq;q^2)_{j+1} (q)_j}
  \label{d2k3i2}
\end{gather}
\end{cor}
  
 \begin{pf}
    To obtain~\eqref{d2k3i1}, replace $a$ by $aq^2$ in~\eqref{d2k3i3}.
To obtain~\eqref{d2k3i2}, subtract $a^2 q^4$ times \eqref{d2k3i1}
with $a$ replaced by $aq^2$ from \eqref{d2k3i3}.
  \qed\end{pf}

Not surprisingly, other limiting cases of~\eqref{gen14} reduce a particular double series to a familiar single sum.
\begin{cor}  
\begin{equation} \label{sl79}
 (-q)_\infty \sum_{j,k\geqq 0} \frac{q^{4j^2+6k^2+8jk-k} }{(-q;q^4)_{j+k} (q^4;q^4)_j (q^4;q^4)_k}
 =\sum_{j=0}^\infty \frac{q^{j^2}}{(q)_{2j}}
 \end{equation}
 \end{cor}
 \begin{pf} In~\eqref{gen14}, let $b,y,n\to\infty$, 
set $x=-\sqrt{q}$, $a=1$, and replace $q$ by $q^2$ throughout.\qed\end{pf}
 
% \begin{rem} 
The right hand side of~\eqref{sl79} appears twice on Slater's 
list~\cite[p. 160, Eq. (79)
 and p. 162, Eq. (98)]{LJS1952}, 
as the series expansion of the (equivalent) infinite products
 $(q^2;q^2)_\infty^{-1} (-q;q^2)_\infty (q^8, q^{12}, q^{20}; q^{20})_\infty$ and
 $(q)_\infty^{-1} (q^2, q^8, q^{10}; q^{10})_\infty (q^6, q^{14}; q^{20})_\infty$ 
 respectively.  This series expansion on the right hand side of~\eqref{sl79} 
 is originally due to L.J. Rogers~\cite[p. 330, 2nd eq.]{LJR1894}.
% \end{rem}

%\begin{rem}
  Andrews~\cite{GEA:PC} pointed out the following alternate
simplification of the double sum in the left hand side of
\eqref{sl79}:
  \begin{eqnarray*}
  &&  \sum_{j,k\geqq 0} \frac{q^{4j^2+6k^2+8jk-k} }
     {(-q;q^4)_{j+k} (q^4;q^4)_j (q^4;q^4)_k} \\
   &=& \sum_{t=0}^\infty \frac{q^{4t^2} }{(-q^4;q^4)_t (q^4;q^4)_t}
    \sum_{k=0}^t \frac{q^{2k^2-k} (q^4;q^4)_t}{(q^4;q^4)_k (q^4;q^4)_{t-k} } 
    \mbox{\qquad (by letting $t=j+k$)} \\
   &=& \sum_{t=0}^\infty \frac{ q^{4t^2}}{(q^4;q^4)_t} 
     \mbox{\qquad (by the $q$-binomial theorem~\cite[p. 8, Eq. (1.3.2)]{GR}). }
  \end{eqnarray*}
Notice that the last expression is the series portion of the first 
Rogers-Ramanujan identity (the $k=i=2$ case of \eqref{AndGor})
with $q\to q^4$.
%\end{rem}

Next, consider the $k=4$ case of Andrews' transformation:
\begin{multline}
\vwp{12}{11}{a}{b,c,d,e,f,g,h,i,q^{-n}}{q}{a^4 q^{n+4}/bcdefghi}\\
= \frac{(aq)_n (aq/hi)_n}{(aq/h)_n (aq/i)_n}
  \sum_{j=0}^n \frac{(q^{-n})_j (h)_j (i)_j (aq/fg)_j q^j}
                {(q)_j (aq/f)_j (aq/g)_j (hiq^{-n}/a)_j} \\ \times
  \sum_{k=0}^j \frac{(q^{-j})_k (f)_k (g)_k (aq/de)_k q^k}
                {(q)_k (aq/d)_k (aq/e)_k (fgq^{-j}/a)_k} 
  \hgs{4}{3}{q^{-k},d,e,aq/bc}{aq/b,aq/c,deq^{-j}/a}{q}{q} \label{GEAk4}
\end{multline}
and
Verma and Jain's transformation~\cite[p. 232, Eq. (1.5)]{VJ}:
\begin{gather}
\vwp{12}{11}{a}{x,xq,xq^2,y,yq,yq^2,q^{2-n},q^{1-n},q^{-n}}{q^3}{a^4 q^{3n+3}/x^3 y^3} \nonumber\\
= \frac{(aq)_n(aq/xy)_n}{(aq/x)_n (aq/y)_n} 
   \hgs{6}{5}{\sqrt[3]{a} , \omega \sqrt[3]{a}, \omega^2 \sqrt[3]{a},x,y,q^{-n}}
{\sqrt{a},-\sqrt{a},\sqrt{aq},-\sqrt{aq},xy q^{-n}/a}{q}{q},
\label{VJ12-11-2}
\end{gather}
where $\omega$ is a primitive cube root of unity.
\begin{thm} Identity~\eqref{gen27} is valid.
\end{thm}
\begin{pf}
The proof is completely analogous to that of identity~\eqref{gen14}, with
\eqref{GEAk4} playing the role of \eqref{GEAk3}, and \eqref{VJ12-11-2}
playing the role of \eqref{VJ10-9-2}.  This time the ``very-well-poised link" 
is
\[ \vwp{12}{11}{a}{x,xq,xq^2,y,yq,yq^2, q^{2-3n}, q^{1-3n}, q^{-3n}}{q^3}{a^4 q^{9n+3}/x^3 y^3}.\]
\qed\end{pf}

\begin{cor} Identity~\eqref{d3k4i4} is valid. \end{cor}
\begin{pf} Let $b, x, y, n \to \infty$ in Eq.~\eqref{gen27}. \qed\end{pf}

Just like~\eqref{d2k3i3}, Eq.~\eqref{d3k4i4} is one of a set of
closely related identities; the three partners of~\eqref{d3k4i3} are
\begin{gather}
\sum_{r=0}^\infty \sum_{s=0}^\infty \sum_{t=0}^\infty
     \frac{a^{r+2s+3t} q^{3[(r+s+t)^2 + (s+t)^2 + t^2 +3r+2s+t]}}
{(q^3;q^3)_r (q^3;q^3)_s (q^3;q^3)_t}
     =\frac{(aq)_\infty}{(aq^3;q^3)_\infty}\sum_{j=0}^\infty 
\frac{a^j q^{j^2+3j} (aq^3;q^3)_j}{(aq)_{2j+2} (q)_j}, \label{d3k4i1}\\
\sum_{r=0}^\infty \sum_{s=0}^\infty \sum_{t=0}^\infty
     \frac{a^{r+2s+3t} q^{3[(r+s+t)^2 + (s+t)^2 + t^2 + 2r+2s+t]}}
{(q^3;q^3)_r (q^3;q^3)_s (q^3;q^3)_t}
     =\frac{(aq)_\infty}{(aq^3;q^3)_\infty}
\sum_{j=0}^\infty \frac{a^j q^{j^2+2j} (aq^3;q^3)_j}{(aq)_{2j+2} (q)_j},
\label{d3k4i2}\\
\sum_{r=0}^\infty \sum_{s=0}^\infty \sum_{t=0}^\infty
     \frac{a^{r+2s+3t} q^{3[(r+s+t)^2 + (s+t)^2 + t^2 + r+s+t]}}
{(q^3;q^3)_r (q^3;q^3)_s (q^3;q^3)_t}
     =\frac{(aq)_\infty}{(aq^3;q^3)_\infty}
\sum_{j=0}^\infty \frac{a^j q^{j^2+j} (aq^3;q^3)_j}{(aq)_{2j+1} (q)_j}.
\label{d3k4i3}
\end{gather}

\section{Additional Results}\label{AddResults}
  In light of the previous section, it makes sense to consider 
another transformation formula of
Verma and Jain~\cite[p. 232, Eq. (1.3)]{VJ}; 
see Bailey~\cite[p. 6, Eq. (6.1)]{WNB1948} for the $n\to\infty$ case.
\begin{gather}\label{VJ10-9-1}
   \vwp{10}{9}{a}{b,x,-x,y,-y,-q^{-n},q^{-n}}{q}{-a^3 q^{2n+3}/bx^2y^2} \\=
   \frac{(a^2 q^2;q^2)_n (a^2 q^2/x^2 y^2;q^2)_n}{(a^2 q^2/x^2;q^2)_n (a^2 q^2/y^2)_n}
       \hgs{5}{4}{ x^2, y^2, -aq/b, -aq^2/b, q^{-2n}}{-aq,-aq^2, a^2 q^2/b^2, x^2 y^2 q^{-2n}/a^2}{q^2}{q^2}
       \nonumber .
 \end{gather}
 Here all that is needed is to make the substitutions $c=x$, $d=-x$, $e=y$, $f=-y$, $g=-q^{-n}$ in \eqref{GEAk3},
 equate its right hand side with the right hand side of~\eqref{VJ10-9-1}, and after some routine algebra,
 results in the following theorem:
 \begin{thm}
 \begin{multline}
 \sum_{j=0}^n \frac{ (q^{-2n};q^2)_j (-y)_j (-aq/xy)_j q^j  }{ (q)_j (-aq/x)_j (aq/y)_j (yq^{-2n}/a)_j  }
  \hgs{4}{3}{q^{-j},-x,y,aq/bx}{aq/b, aq/x, -xyq^{-j}/a}{q}{q}\\
  = \frac{ (-aq)_{2n} (a^2 q^2/ x^2 y^2;q^2)_n }{ (a^2 q^2/x^2;q^2)_n (aq/y)_{2n} }
  \hgs{5}{4}{q^{-2n},x^2,y^2, -aq/b, -aq^2/b}{x^2 y^2 q^{-2n}/a^2, a^2 q^2/b^2, -aq, -aq^2}{q^2}{q^2}
 , \end{multline}
 \end{thm} 
\noindent which, after suitable specialization, yields  
 \begin{cor}
\begin{equation}\label{aRS}
    \sum_{j,k\geqq 0} \frac{ (-1)^k a^{j+2k} q^{(j+k)^2 + \binomial{k+1}{2}}}{(q)_j (q)_k} 
    = (-aq)_\infty \sum_{j=0}^\infty 
\frac{ a^{2j} q^{2j^2}}{(-aq)_{2j} (q^2;q^2)_j}.
 \end{equation}
 \end{cor}
 
% \begin{rem} 
The series on the right hand side of~\eqref{aRS} 
with $a=1$ is the series associated with
 the first Rogers-Selberg identity, an expansion of the series 
 $(q^3,q^4,q^7;q^7)_\infty (q^2;q^2)_\infty^{-1}$, due to Rogers~\cite[p. 338]{LJR1894}
 and recorded by Slater~\cite[Eq. (33)]{LJS1952}.
% \end{rem}

\section{Discussion}\label{Disc}
  Once I had in hand a $q$-hypergeometric explanation for the existence
of identities like~\eqref{d2k3i3} and~\eqref{d3k4i4}, it was only
natural to look for additional analogous identities.  I did not 
search exhaustively, but rather presented a couple of striking examples
relating to well-known series (e.g. Rogers-Selberg).  Certainly additional
identities of this type exist (e.g. a multisum version of the Bailey
``mod 9 identities"~\cite[p. 422, Eqs. (1.6)--(1.8)]{WNB1947},
~\cite[Eqs. (40)--(42)]{LJS1952}), and the 
interested reader is encouraged to use the methods of this paper to
work out additional examples. 

  A more ambitious project would be to look for bijective proofs of
identities like~\eqref{d2k3i3} and \eqref{d3k4i4}. 
 
\section*{Acknowledgements}
The author thanks the anonymous referee and the editors,
Joseph Kung and Doron Zeilberger, for carefully reading the
manuscript and providing a number of helpful suggestions.

\end{document}